\documentclass[11pt,a4paper]{article}

\usepackage[T1]{fontenc}
\usepackage[utf8]{inputenc}
\usepackage{lmodern}
\usepackage[margin=28mm]{geometry}
\usepackage{microtype}

\usepackage{amsmath,amssymb,amsthm}
\usepackage{booktabs}
\usepackage{enumitem}
\usepackage{placeins}

\usepackage[hidelinks]{hyperref}
\hypersetup{
  pdftitle={Weakly intersecting--union families},
  pdfauthor={},
  pdfsubject={Extremal set theory}
}

\numberwithin{equation}{section}
\theoremstyle{plain}
\newtheorem{theorem}{Theorem}[section]
\newtheorem{lemma}[theorem]{Lemma}

\theoremstyle{definition}

\newcommand{\m}{\mathcal }

\newcommand{\F}{\mathcal F}
\newcommand{\I}{\mathcal I}
\newcommand{\D}{\mathcal D}
\newcommand{\Z}{\mathcal Z}
\renewcommand{\H}{\mathcal H}

\newcommand{\one}{\mathbf 1}
\DeclareMathOperator{\rank}{rank}
\DeclareMathOperator{\col}{col}

\setlist[enumerate]{label=\textup{(\roman*)},leftmargin=*,itemsep=3pt}
\allowdisplaybreaks

\begin{document}

\title{Weakly intersecting--union families}
\author{
Yang Huang\thanks{Moscow Institute of Physics and Technology; 
E-mail: \url{1060393815@qq.com}.}
\quad \quad 
Andrey Kupavskii\thanks{Moscow Institute of Physics and Technology; 
E-mail: \url{kupavskii@ya.ru}. }
\quad \quad
Anyuan Tian\thanks{Moscow Institute of Physics and Technology; 
E-mail: \url{tyan.a@mipt.ru}. }
}
\date{}

\maketitle

\vspace{-0.5cm}

\begin{abstract}
A family of subsets of an $n$-element set, $n\ge2$, is weakly intersecting--union if any two members intersect and fail to cover the ground set unless they are complementary. We determine its maximum size, resolving the last remaining problem in Gronau's 1981 classification of Boolean conditions on pairs of sets. The maximum is $2^{n-2}$ for odd $n$ and for $n\ge10$, and $\binom n{n/2}$ for $n=2,4,6,8$. In the proof, we use the Marica--Sch\"onheim difference inequality, bounds on ranks of incidence matrices, and a double count of common subsets that omit a given point.

\end{abstract}
\section{Introduction}\label{sec:intro}
Write $[n]=\{1,\ldots,n\}$, and assume $n\ge2$. For a set $X$, let $2^X$ denote its power set. A family $\F\subseteq2^{[n]}$ is \emph{intersecting--union}, or \emph{IU}, if
\[
F\cap G\ne\varnothing\quad\text{and}\quad F\cup G\ne[n]\qquad(F,G\in\F).
\]
It is well-known that a size of such $\F$ is at most $2^{n-2}$. One short proof uses the Harris--Kleitman correlation inequality~\cite{Harris,Kleitman}. Let $\mathcal A$ be the upward closure of $\F$ and $\mathcal B$ its downward closure. The family $\mathcal A$ is intersecting, and no two members of $\mathcal B$ cover $[n]$. Each therefore contains at most one member of every complementary pair, giving $|\mathcal A|,|\mathcal B|\le2^{n-1}$. Since an increasing and a decreasing family are negatively correlated under the uniform measure,
\[
|\F|\le|\mathcal A\cap\mathcal B|\le\frac{|\mathcal A||\mathcal B|}{2^n}\le2^{n-2}.
\]
The family $\{F\subseteq[n]:1\in F,\ n\notin F\}$ attains this bound.

Gronau~\cite{Gronau} suggested to study the following weaker property. A family is \emph{weakly intersecting--union}, or \emph{weakly IU}, if
\[
F\cap G\ne\varnothing\quad\text{and}\quad F\cup G\ne[n]\qquad(F,G\in\F,\ G\ne[n]\setminus F).
\]
That is, the property is the same as IU, except we additionally allow complementary sets.

Of course, any IU family is also weakly IU, and the example above provides us with a weakly IU family of size $2^{n-2}$. For small even $n$, we can do better: it is easy to see that the middle level $\binom{[n]}{n/2}$ of the hypercube is weakly IU. Its size exceeds $2^{n-2}$ exactly for $n\in\{2,4,6,8\}$.

Gronau conjectured that one of the two constructions must be extremal (cf. also Frankl~\cite[Conjecture~11.5]{Frankl}). In fact, Gronau~\cite{Gronau} considered extremal families defined by Boolean conditions on pairs of sets, and the maximum size of a weakly IU family was the sole undetermined maximum in Gronau's classification; see also Gronau--Proske~\cite[p.~441, footnote~4]{GP}. We prove this conjecture.
\begin{theorem}\label{thm:main}
Let $n\ge2$, and let $\F\subseteq2^{[n]}$ be weakly IU. Then
\[
|\F|\le\begin{cases} \binom n{n/2},&n\in\{2,4,6,8\},\\
2^{n-2},&n\text{ is odd or }n\ge10.
\end{cases}
\]
\end{theorem}

\subsection{Proof sketch}\label{sec:sketch}
We first show that shifting preserves the weak IU property and assume that an extremal weakly IU family $\F$ is shifted. We then construct a decomposition of $\F$ into four families $\F_{00},\F_{01},\F_{10},\F_{11}$ on $R=\{2,\ldots,n-1\}$ according to the membership of the elements $1$ and $n$ in a set. For a family $\mathcal U$ on $R$, let $\mathcal U^c$ consist of the complements of its members in $R$. 

The first key structural step is to show  that the four families $\F_{10}^c$, $\F_{00}\cup\F_{11}$, $\D$, and $\D^c$ are disjoint. Here $\D$ consists of the sets that are strict subsets of both a member of $\F_{11}$ and a member of $\F_{00}^c$. The proof uses shiftedness and different structural properties of weakly IU families. 

An easy consequence of the previous step is the following inequality.
\[
|\F|\le2^{n-2}+|\F_{01}|+|\I|-2|\D|,
\]
where $\I = \F_{00}\cap \F_{11}$. We then use the Marica--Sch\"onheim inequality  \cite{MS} 
$|\{P\setminus Q:P,Q\in\H\}|\ge|\H|$  and some additional structural and linear-algebraic information to show that $2|\D|\ge |\F_{01}|+|\I|$ (except when $n=2,4,6,8$).

In fact, in the case when $n$ is odd, there is a decomposition of sets from $\F_{01}$ and $\I$ according to their sizes into two families $\m H_{0}$ and $\m H_{1}$, such that $\{P\setminus Q:P,Q\in\H_i\}\subseteq\D$  holds for $i\in \{0,1\}$, and then the Marica--Sch\"onheim inequality finishes the proof.

When $n$ is even, the problematic sets are sets from $\m F_{01}$ of size $k=(n-2)/2$: they cannot be easily included in $\m H_i$. We then develop a linear-algebraic argument that consists of two steps. First, we prove a two-family refinement of the  Marica--Sch\"onheim inequality. It gives   $|\D|\ge|\H_i|+\rank M$ for a certain inclusion matrix $M$. The second step, also partly linear-algebraic, is to relate $\rank M$ to the number of problematic sets from the middle layer.

The cases of small even $n$ use the same argument but require a bit more care.

Section~\ref{sec:decomposition} gives the decomposition and the proof for odd $n$. Section~\ref{sec:central} establishes the central family's combinatorial properties. Section~\ref{sec:linear} proves the linear algebra lemmas. Section~\ref{sec:even} completes the proof of Theorem~\ref{thm:main} in all even dimensions.

\section{Decomposition and the odd case}\label{sec:decomposition}
The case $n=2$ of Theorem~\ref{thm:main} is immediate, since only the two singletons may belong to $\F$. Henceforth assume $n\ge3$.

For $x\notin S$ and $y\in S$, write $S+x=S\cup\{x\}$ and $S-y=S\setminus\{y\}$. For distinct $i,j\in[n]$, define the $(i,j)$-shift of $\F$ by
\[
S_{ij}(\F)=\{S_{ij}(F):F\in\F\},\qquad S_{ij}(F)=\begin{cases} F-j+i,&j\in F,\ i\notin F,\ F-j+i\notin\F,\\
F,&\text{otherwise}.
\end{cases}
\]
A family is \emph{shifted} if $S_{ij}(\F)=\F$ for every $i<j$.
\begin{lemma}\label{lem:shift}
Shifting preserves size and the weakly IU property. Consequently, every weakly IU family can be replaced by a shifted weakly IU family of the same size.
\end{lemma}
\begin{proof}
Size preservation is immediate. Suppose that $F,G\in S_{ij}(\F)$ are disjoint and $|F|+|G|<n$. At least one, say $F$, is not in $\F$, so $i\in F$, $j\notin F$ and $F'=F-i+j\in\F$. Since $i\notin G$, we have $G\in\F$. If $j\notin G$, then $F'$ and $G$ are disjoint. If $j\in G$, the fact that $G$ was not shifted implies $G'=G-j+i\in\F$, and $F'$ and $G'$ are disjoint. In either case we obtain two members of $\F$ whose total size is less than $n$, contradicting the weakly IU property. Taking complements in $[n]$ interchanges $(i,j)$-shifts and $(j,i)$-shifts, so the same argument rules out non-complementary covering pairs. For an IU family the disjoint-pair argument also applies when $|F|+|G|=n$, since the resulting pair in $\F$ is still disjoint. Taking complements then excludes all covering pairs.

Every nontrivial shift with $i<j$ decreases $\sum_{F\in\F}\sum_{x\in F}x$. Repeatedly applying such shifts therefore yields a shifted family.
\end{proof}
We now work with a shifted family $\F$. Put $R=\{2,\ldots,n-1\}$. Complements of subsets of $R$ are taken in $R$; thus $\overline S=R\setminus S$ and $\mathcal U^c=\{\overline S:S\in\mathcal U\}$ for $\mathcal U\subseteq2^R$. Define the four sections, with the two indices recording membership of $1$ and $n$, respectively,
\begin{equation}\label{eq:sections}
\begin{aligned}
\F_{10}&=\{S\subseteq R:S+1\in\F\},&\F_{00}&=\{S\subseteq R:S\in\F\},\\
\F_{01}&=\{S\subseteq R:S+n\in\F\},&\F_{11}&=\{S\subseteq R:S+1+n\in\F\}.
\end{aligned}
\end{equation}
Set
\begin{equation}\label{eq:D}
\I=\F_{00}\cap\F_{11},\qquad \D=\{S\subseteq R:S\subsetneq D,\ S\subsetneq\overline B\text{ for some }D\in\F_{11},\ B\in\F_{00}\}.
\end{equation}
The family $\D$ is a downset: if $S\in\D$ and $T\subseteq S$, the same two witnesses in \eqref{eq:D} give $T\in\D$.
\begin{lemma}\label{lem:four}
The families in \eqref{eq:sections}--\eqref{eq:D} have the following properties.
\begin{enumerate}[label=\textup{(\roman*)}]
\item\label{item:four-comparable} Comparable members of $\F_{11}$ and $\F_{00}^c$ are equal.
\item The four families $\F_{10}^c$, $\F_{00}\cup\F_{11}$, $\D$, and $\D^c$ are pairwise disjoint.
\item \label{item:nondifference}If $X,Y\in\F_{01}\cup\I$ and $X\setminus Y\notin\D$, then
\[
X,Y\in\F_{01}\setminus\I,\qquad |X|+|Y|=n-2,\qquad\overline X,\overline Y\in\F_{10}.
\]
\end{enumerate}
\end{lemma}
\begin{proof}
The family $\F_{00}$ is intersecting, since its members omit $1,n$. Similarly, no two members of $\F_{11}$ 
cover $R$: the corresponding members $D+1+n$ of $\F$ share $1,n$. Moreover, the weakly IU property gives
\begin{equation}\label{eq:mixed}
\begin{aligned}
X\cap B&\ne\varnothing&& (X\in\F_{10}\cup\F_{01},\ B\in\F_{00}),\\
X\cup D&\ne R&&(X\in\F_{10}\cup\F_{01},\ D\in\F_{11}).
\end{aligned}
\end{equation}
For the first assertion, choose $e\in\{1,n\}$ with $X+e\in\F$. The sets $X+e$ and $B$ omit the element of $\{1,n\}$ distinct from $e$, so they are non-complementary and must intersect. For the second, $X+e$ and $D+1+n$ share $e$, so the weakly IU property prevents their union from being $[n]$.

\textup{(i)}. If $D\in\F_{11}$ and $B\in\F_{00}$ satisfy $D\subsetneq\overline B$, then $D+1+n$ and $B$ are disjoint non-complementary members of $\F$. If $\overline B\subsetneq D$, the sets $D+1+n$ and $B$ cover $[n]$ without being complementary. Both contradict the weakly IU property.

\textup{(ii)}. By \eqref{eq:mixed}, no member of $\F_{10}^c$ is contained in a member of $\F_{11}$ or contains a member of $\F_{00}$. Every member of $\F_{11}\cup\D$ is contained in some member of $\F_{11}$; every member of $\F_{00}\cup\D^c$ contains some member of $\F_{00}$. Thus $\F_{10}^c$ is disjoint from the other three families.

For $S\in\D$, choose $D\in\F_{11}$ and $B\in\F_{00}$ with $S\subsetneq D$ and $S\subsetneq\overline B$. If $S\in\F_{11}$, the latter containment contradicts (i); if $\overline S\in\F_{00}$, the former does. Thus $\D\cap\F_{11}=\varnothing$ and $\D^c\cap\F_{00}=\varnothing$. Since $\F_{00}$ is intersecting and every member of $\D$ is disjoint from some member of $\F_{00}$, we have $\D\cap\F_{00}=\varnothing$. Since no two members of $\F_{11}$ cover $R$ and every member of $\D$ is contained in a member of $\F_{11}$, we have $\D^c\cap\F_{11}=\varnothing$ and $\D\cap\D^c=\varnothing$.

\textup{(iii)}. Shifting gives
\begin{equation}\label{eq:shift}
\F_{01}\subseteq\F_{10},\quad C-c\in\F_{11}\ (c\in C\in\F_{01}),\quad C+e\in\F_{00}\ (C\in\F_{01},\ e\in R\setminus C).
\end{equation}
These follow by replacing, respectively, $n$ by $1$, $c$ by $1$, and $n$ by $e$ in $C+n$.

If $X,Y\in\F_{01}$, apply the weakly IU property to $X+1,Y+n$, using $\F_{01}\subseteq\F_{10}$. This shows that $\F_{01}$ is weakly IU on $R$. If $X\in\F_{01}$ and $Y\in\I$, then $Y\in\F_{00}\cap\F_{11}$, so \eqref{eq:mixed} gives $X\cap Y\ne\varnothing$ and $X\cup Y\ne R$. If $X,Y\in\I$, the same conclusions follow because $\F_{00}$ is intersecting and no two members of $\F_{11}$ cover $R$. These cases show that $\F_{01}\cup\I$ is weakly IU on $R$, and that any complementary pair $X,Y$ in this family has both members in $\F_{01}\setminus\I$ and satisfies $|X|+|Y|=n-2$. For such a pair, we also have $\overline X=Y\in\F_{01}\subseteq\F_{10}$, and similarly $\overline Y\in\F_{10}$.

It remains to consider a non-complementary pair $X,Y\in\F_{01}\cup\I$. Put $S=X\setminus Y$ and assume $S\notin \D$. Choose $c\in X\cap Y$ and $e\in R\setminus(X\cup Y)$, and define
\[
D_X=\begin{cases}X,&X\in\I,\\X-c,&X\notin\I,\end{cases} \qquad B_Y=\begin{cases}Y,&Y\in\I,\\Y+e,&Y\notin\I.\end{cases}
\]
By the definition of $\I$ and \eqref{eq:shift}, we have $D_X\in\F_{11}$ and $B_Y\in\F_{00}$. The choices of $c$ and $e$ also give
\[
S\subseteq D_X,\qquad S\subseteq\overline{B_Y}.
\]
Since $S\notin\D$, at least one of these containments is an equality. Thus $D_X$ and $\overline{B_Y}$ are comparable, so (i) gives
\[
D_X=S=\overline{B_Y}.
\]
These equalities force $X,Y\notin\I$: if $X\in\I$, then $D_X=X$, contradicting $c\in X\setminus S$; if $Y\in\I$, then $\overline{B_Y}=\overline Y$, contradicting $e\in\overline Y\setminus S$. Hence $X,Y\in\F_{01}\setminus\I$, and
\[
X=D_X+c=S+c,\qquad\overline Y=\overline{B_Y}+e=S+e,\qquad |X|-1=|S|=(n-2)-|Y|-1.
\]
Finally, we check that $\overline X,\overline Y\in\F_{10}$. We have $\overline S=B_Y\in\F_{00}$ and $S=D_X\in\F_{11}$. Shifting $c$ to $1$ in $\overline S$ gives $\overline S-c=\overline X\in\F_{10}$; shifting $n$ to $e$ in $S+1+n$ gives $S+e=\overline Y\in\F_{10}$. This proves (iii).
\end{proof}
The four sections in \eqref{eq:sections} give the first equality below. The second uses $\I=\F_{00}\cap\F_{11}$. For the inequality, apply Lemma~\ref{lem:four}(ii): its four disjoint families are contained in $2^R$, so $|\F_{10}|+|\F_{00}\cup\F_{11}|+2|\D|\le2^{n-2}$. Thus
\begin{equation}\label{eq:count}
\begin{aligned}
|\F|&=|\F_{00}|+|\F_{01}|+|\F_{10}|+|\F_{11}|\\
&=|\F_{10}|+|\F_{00}\cup\F_{11}|+|\F_{01}|+|\I|\\
&\le2^{n-2}+|\F_{01}|+|\I|-2|\D|.
\end{aligned}
\end{equation}

For odd $n$ and for $n\geq 10$,  it suffices to prove $|\F_{01}|+|\I|\le2|\D|$ to settle the main part of Theorem~\ref{thm:main}. In the exceptional cases $n=4,6,8$ we will instead prove $|\F_{01}|+|\I|-2|\D| \le\binom n{n/2}-2^{n-2}.$
\begin{proof}[Proof of Theorem~\ref{thm:main} for odd $n$]
Since $n-2$ is odd, for $\mathcal U\in\{\F_{01},\I\}$ we have the partition
\[
\mathcal U^-=\{U\in\mathcal U:|U|<(n-2)/2\},\qquad \mathcal U^+=\{U\in\mathcal U:|U|>(n-2)/2\}.
\]
Put
\[
\H_0=\F_{01}^-\cup\I^+,\qquad\H_1=\F_{01}^+\cup\I^-.
\]
We claim that $P\setminus Q\in\D$ for every $P,Q\in\H_0$ and for every $P,Q\in\H_1$. Suppose otherwise. By Lemma~\ref{lem:four}(iii), $P,Q\in\F_{01}\setminus\I$ and $|P|+|Q|=n-2$. If $P,Q\in\H_0$, then $P,Q\in\F_{01}^-$, so $|P|+|Q|<n-2$. If $P,Q\in\H_1$, then $P,Q\in\F_{01}^+$, so $|P|+|Q|>n-2$. Both cases contradict the required equality. Applying the Marica--Sch\"onheim inequality to $\H_0$ and $\H_1$ gives
\[
|\H_0|\le|\D|,\qquad|\H_1|\le|\D|.
\]
The two parts of each family lie in opposite size ranges. Hence $ |\F_{01}|+|\I|=|\H_0|+|\H_1|\le2|\D|. $ \eqref{eq:count} now gives $|\F|\le2^{n-2}$, completing the odd case.
\end{proof}

\section{The central level}\label{sec:central}
We now assume that $n$ is even and write $n-2=2k$. If $X,Y\in\F_{01}$, $|X|,|Y|\le k$ and $X\setminus Y\notin\D$, then Lemma~\ref{lem:four}(iii) gives $|X|+|Y|=2k$ and $\overline X,\overline Y\in\F_{10}$. It follows that $|X|=|Y|=k$. Motivated by this, define
\begin{equation}\label{eq:Z}
\Z=\{Z\in\F_{01}:|Z|=k,\ \overline Z\in\F_{10}\}.
\end{equation}
\begin{lemma}\label{lem:central}
The following statements hold.
\begin{enumerate}[label=\textup{(\roman*)}]
\item\label{item:central-incomparable} Every $Z\in\Z$ is incomparable with every distinct member of $\F_{01}\cup\I$. Moreover, $\Z\cap\I=\varnothing$.
\item\label{item:central-cross} If $P\in(\F_{01}\setminus\Z)\cup\I$ and $Z\in\Z$, then
\[
P\setminus Z,\ Z\setminus P\in\D,\qquad |P\cap Z|\ge2,\qquad |R\setminus(P\cup Z)|\ge2.
\]
\end{enumerate}
\end{lemma}
The bounds by $2$ in (ii) use both shifting and the membership of $P\setminus Z$ and $Z\setminus P$ in $\D$. If $|P\cap Z|=1$ or $|R\setminus(P\cup Z)|=1$, one of the differences $P\setminus Z$ and $Z\setminus P$ is a one-point deletion of $Z$ or the complement of a one-point extension of $Z$. Shifting puts that set in $\F_{11}$ or $\F_{00}^c$, respectively, whereas $\D$ is disjoint from both families. These restrictions will be used in Section~\ref{sec:even} when a member of size different from $k$ is present.
\begin{proof}
\textup{(i)}. Let $Z\in\Z$. Since $\overline Z\in\F_{10}$, \eqref{eq:mixed} shows that $Z$ contains no member of $\F_{00}$ and is contained in no member of $\F_{11}$. This proves $Z\notin\I$ and incomparability with members of $\I$. Next, take $P\in\F_{01}$. By \eqref{eq:Z}, $\overline Z\in\F_{10}$, so \eqref{eq:sections} gives $\overline Z+1,P+n\in\F$. If $P\subsetneq Z$, they are disjoint and non-complementary; if $Z\subsetneq P$, they cover $[n]$ and are non-complementary. Both possibilities contradict the weakly IU property.

\textup{(ii)}. Fix $P\in(\F_{01}\setminus\Z)\cup\I$ and $Z\in\Z$. Part (i) gives $\Z\cap\I=\varnothing$, so $P\notin\Z$. If $P\setminus Z\notin\D$, applying Lemma~\ref{lem:four}(iii) with $(X,Y)=(P,Z)$ gives $P\in\F_{01}\setminus\I$, $|P|=2k-|Z|=k$, and $\overline P\in\F_{10}$, which implies $P\in\Z$, a contradiction. If $Z\setminus P\notin\D$, the same lemma applied with $(X,Y)=(Z,P)$ gives the same three conclusions about $P$, again a contradiction. Thus $P\setminus Z,Z\setminus P\in\D$.

By \eqref{eq:Z} and \eqref{eq:shift}, both $Z$ and $\overline Z$ belong to $\F_{10}^c$, so neither belongs to $\D$ by Lemma~\ref{lem:four}(ii). Since $Z\setminus P,P\setminus Z\in\D$, it follows that $P\cap Z\ne\varnothing$ and $P\cup Z\ne R$.

If $P\cap Z=\{z\}$, then $Z\setminus P=Z-z\in\F_{11}$ by \eqref{eq:shift}. This contradicts $Z\setminus P\in\D$, because $\D\cap\F_{11}=\varnothing$ by Lemma~\ref{lem:four}(ii). If $R\setminus(P\cup Z)=\{e\}$, then $P\setminus Z=\overline{Z+e}\in\F_{00}^c$ by \eqref{eq:shift}. This contradicts $P\setminus Z\in\D$, because $\D\cap\F_{00}^c=\varnothing$ is equivalent to $\D^c\cap\F_{00}=\varnothing$, also given by Lemma~\ref{lem:four}(ii). It follows that $|P\cap Z|\ge2$ and $|R\setminus(P\cup Z)|\ge2$.
\end{proof}
For a family $\mathcal U\subseteq2^R$, put
\[
\mathcal U^{<k}=\{U\in\mathcal U:|U|<k\},\qquad\mathcal U^{\le k}=\{U\in\mathcal U:|U|\le k\},
\]
and define $\mathcal U^k$ and $\mathcal U^{>k}$ by the conditions $|U|=k$ and $|U|>k$, respectively. Set
\begin{equation}\label{eq:H}
\begin{aligned}
\H_0&=\F_{01}^{<k}\cup(\F_{01}^k\setminus\Z)\cup\I^{>k},\\
\H_1&=\F_{01}^{>k}\cup\I^{\le k}.
\end{aligned}
\end{equation}
The parts in each union lie in different size ranges, and hence
\begin{equation}\label{eq:Hcount}
|\H_0|+|\H_1|=|\F_{01}|+|\I|-|\Z|.
\end{equation}
We have
\begin{equation}\label{eq:diff}
P\setminus Q\in\D\qquad(P,Q\in\H_i,\ i\in\{0,1\}).
\end{equation}
To prove this, suppose $P,Q\in\H_i$ and $P\setminus Q\notin\D$. Lemma~\ref{lem:four}(iii) gives $P,Q\in\F_{01}\setminus\I$, $|P|+|Q|=2k$, and $\overline P,\overline Q\in\F_{10}$. If $i=1$, then $|P|,|Q|>k$, contradicting $|P|+|Q|=2k$. If $i=0$, then $|P|,|Q|\le k$, so their sum is $2k$ only if $|P|=|Q|=k$. Together with $\overline P,\overline Q\in\F_{10}$, this implies $P,Q\in\Z$, contrary to the definition of $\H_0$.

For a family $\mathcal U\subseteq2^R$, write $\partial^{\le r}\mathcal U$ for the subsets of its members having size at most $r$. The following lemma establishes a property of $\D$ needed for our later application of Lemma~\ref{lem:recover}.
\begin{lemma}\label{lem:small}
Suppose $k\ge2$, and put $\mathcal A=(\F_{01}\cup\I)^k$. Then
\begin{equation}\label{eq:common}
(\partial^{\le k-2}\mathcal A)\cap (\partial^{\le k-2}\mathcal A^c)\subseteq\D.
\end{equation}
\end{lemma}
\begin{proof}
Take $S\subseteq A\setminus B$ with $A,B\in\mathcal A$ and $|S|\le k-2$. Choose $a\in A\setminus S$ and $b\in R\setminus(B\cup S)$, and put
\[
D_0=\begin{cases}A,&A\in\I,\\A-a,&A\notin\I,\end{cases} \qquad B_0=\begin{cases}B,&B\in\I,\\B+b,&B\notin\I.\end{cases}
\]
If $A\notin\I$, then $A\in\F_{01}$ and $A-a\in\F_{11}$ by \eqref{eq:shift}; otherwise $A\in\I\subseteq\F_{11}$. Thus $D_0\in\F_{11}$. The same two alternatives for $B$, using \eqref{eq:shift} and $\I\subseteq\F_{00}$, give $B_0\in\F_{00}$. Therefore, we have $S\subseteq D_0$ and $S\subseteq\overline{B_0}$. Both $D_0$ and $\overline{B_0}$ have size at least $k-1>|S|$, so both containments are strict. Hence $S\in\D$ by \eqref{eq:D}.
\end{proof}
Since $\Z\subseteq(\F_{01}\cup\I)^k$, this gives in particular
\begin{equation}\label{eq:small}
Z,Z'\in\Z,\quad S\subseteq Z\setminus Z',\quad |S|\le k-2\quad\Longrightarrow\quad S\in\D.
\end{equation}
Taking $Z=Z'$ gives $\varnothing\in\D$ whenever $\Z\ne\varnothing$ and $k\ge2$.

\section{Linear algebra}\label{sec:linear}
This section collects the linear-algebraic results used in Section~\ref{sec:even}. All vector spaces and ranks are over $\mathbb R$. For a finite set $V$ and families $\mathcal A,\mathcal B\subseteq2^V$, write
\[
W(\mathcal A,\mathcal B)=(\one_{B\subseteq A})_{A\in\mathcal A,B\in\mathcal B}.
\]
Thus it is the inclusion matrix with rows indexed by $\mathcal A$ and columns indexed by $\mathcal B$. For $0\le r\le|V|$, we abbreviate $W_r(\mathcal A)=W(\mathcal A,\binom Vr)$. For a matrix $A$ with row set $\Omega$ and $\mathcal U\subseteq\Omega$, write $A[\mathcal U]$ for its restriction to the rows in $\mathcal U$. Both of the next two proofs use the elementary identity
\begin{equation}\label{eq:binomial}
\sum_{T\subseteq E}(-1)^{|T|}=\begin{cases}1,&E=\varnothing,\\0,&E\ne\varnothing,\end{cases}
\end{equation}
which follows from the binomial theorem.

If $P\setminus Q\in\D$ for all $P,Q\in\H$, then the Marica--Sch\"onheim gives $|\D|\ge|\H|$. For a downset $\D$, the following lemma provides a sharpening of this inequality.
\begin{lemma}\label{lem:dimension}
Let $V$ be a finite set, let $\D\subseteq2^V$ be a downset, and let $\H,\mathcal C\subseteq2^V$. Suppose that
\[
P\setminus Q\in\D\quad(P,Q\in\H),\qquad C\setminus Q\in\D\text{ and }C\nsubseteq Q\quad(C\in\mathcal C,Q\in\H).
\]
Then
\[
|\D|\ge|\H|+\rank W(\mathcal C,\D).
\]
\end{lemma}
\begin{proof}
If $\H=\varnothing$, the conclusion says that the rank of $W(\mathcal C,\D)$ is at most its number of columns. Assume $h=|\H|>0$, and define the $|\D|\times h$ matrix
\[
N=\bigl((-1)^{|S|}\one_{S\cap Q=\varnothing}\bigr)_{S\in\D,Q\in\H}.
\]
For $A\in\H\cup\mathcal C$ and $Q\in\H$, we have $A\setminus Q\in\D$ by the two conditions in the statement of the lemma. Since $\D$ is a downset, it contains every subset of $A\setminus Q$. The $(A,Q)$ entry of the matrix product is therefore
\[
\bigl(W(\H\cup\mathcal C,\D)N\bigr)_{A,Q} =\sum_{S\in\D}\one_{S\subseteq A}(-1)^{|S|}\one_{S\cap Q=\varnothing} =\sum_{S\subseteq A\setminus Q}(-1)^{|S|}=\one_{A\subseteq Q},
\]
where the last equality is \eqref{eq:binomial}. Taking the rows in $\H$ gives
\[
W(\H,\D)N=W(\H,\H)^T.
\]
Enumerate $\H=\{H_1,\ldots,H_h\}$ so that $|H_1|\le\cdots\le|H_h|$, using the same order for the rows and columns of $W(\H,\H)^T$. Then $W(\H,\H)^T=(\one_{H_i\subseteq H_j})_{i,j=1}^h$ is upper triangular with diagonal entries $1$, since $H_i\nsubseteq H_j$ whenever $i>j$. It is invertible, and hence $N$ has rank $h$.

For $C\in\mathcal C$ and $Q\in\H$, the same entry calculation gives
\[
\bigl(W(\mathcal C,\D)N\bigr)_{C,Q}=\one_{C\subseteq Q}=0,
\]
where the last equality follows from the assumption $C\nsubseteq Q$. Thus $W(\mathcal C,\D)N=0$. Every transposed row of $W(\mathcal C,\D)$ belongs to $\ker N^T$, so rank--nullity gives
\[
\rank W(\mathcal C,\D)\le\dim\ker N^T=|\D|-\rank N=|\D|-|\H|.\qedhere
\]
\end{proof}

\begin{lemma}\label{lem:recover}
Let $V$ be a finite set, let $r\ge0$ be an integer, and let $\varnothing\ne\mathcal U\subseteq2^V$ and $\D\subseteq2^V$. Suppose that
\begin{equation}\label{eq:recover-condition}
U,U'\in\mathcal U,\quad S\subseteq U\setminus U', \quad |S|\le r\quad\Longrightarrow\quad S\in\D.
\end{equation}
For $S\subseteq V$, put $c_S=(\one_{S\subseteq U})_{U\in\mathcal U}$. Then $c_S\in\col W(\mathcal U,\D)$ whenever $|S|\le r$.
\end{lemma}
\begin{proof}
We proceed by induction on $|S|$. Taking $U'=U$ in \eqref{eq:recover-condition} gives $\varnothing\in\D$, so $c_{\varnothing}$ is a column of $W(\mathcal U,\D)$.

Now let $1\le|S|\le r$ and assume the conclusion holds for all proper subsets of $S$. If $c_S=0$ or $S\in\D$, the conclusion is immediate. Otherwise, choose $U_0\in\mathcal U$ with $S\subseteq U_0$. We claim that $S$ intersects every member of $\mathcal U$. Otherwise, there exists $U\in\mathcal U$ with $S\cap U=\varnothing$. Then $S\subseteq U_0\setminus U$, so \eqref{eq:recover-condition} gives $S\in\D$, a contradiction.

Thus, for every $U\in\mathcal U$, \eqref{eq:binomial} gives
\[
\sum_{T\subseteq S}(-1)^{|T|}\one_{T\subseteq U} =\sum_{T\subseteq S\cap U}(-1)^{|T|}=0.
\]
Isolating the term $T=S$ in each coordinate yields
\begin{equation}\label{eq:recovery}
c_S=(-1)^{|S|+1}\sum_{T\subsetneq S}(-1)^{|T|}c_T.
\end{equation}
Each $c_T$ on the right belongs to $\col W(\mathcal U,\D)$ by induction, and hence so does $c_S$.
\end{proof}

\begin{lemma}\label{lem:restricted-rank}
Let $|V|=v$ and let $p,r$ be integers with $0\le r\le\min(p,v-p)$. For every $\mathcal U\subseteq\binom Vp$,
\begin{equation}\label{eq:restricted-rank}
\rank W_r(\mathcal U)\ge \frac{\binom vr}{\binom vp}|\mathcal U|.
\end{equation}
\end{lemma}
\begin{proof}
By Gottlieb's theorem \cite[Corollary~2]{Gottlieb}, see also \cite[Theorem~23.7]{FT}, the incidence matrix $W_r(\binom Vp)$ has full column rank. Hence we may choose $\mathcal B\subseteq\binom Vp$ with $|\mathcal B|=\binom vr$ such that the corresponding rows are linearly independent. For every permutation $\sigma$ of $V$, let $\sigma(\mathcal B)=\{\sigma(B):B\in\mathcal B\}$. The rows indexed by $\sigma(\mathcal B)$ remain linearly independent, since $\sigma$ merely permutes the rows and columns of $W_r(\binom Vp)$. Consequently, \begin{equation}\label{eq:permutation-rank} \rank W_r(\mathcal U) \ge |\mathcal U\cap\sigma(\mathcal B)|.
\end{equation}

Averaging over all permutations $\sigma$ of $V$, we obtain
\[
\rank W_r(\mathcal U)
\ge \frac1{v!}\sum_\sigma
    |\mathcal U\cap\sigma(\mathcal B)|
=\frac{\binom vr}{\binom vp}|\mathcal U|.\qedhere
\]
\end{proof}

\section{The proof in even dimensions}\label{sec:even}
The case $n=2$ was settled in Section~\ref{sec:decomposition}. Let $n=2k+2$ with $k\ge1$, so $|R|=2k$, and set $M=W(\Z,\D)$. By \eqref{eq:Hcount}, we have $|\F_{01}|+|\I|=|\H_0|+|\H_1|+|\Z|$.

Fix $i\in\{0,1\}$. We verify the conditions of Lemma~\ref{lem:dimension} with $\H=\H_i$ and $\mathcal C=\Z$. First, \eqref{eq:diff} gives $P\setminus Q\in\D$ for all $P,Q\in\H_i$. Next, every $Q\in\H_i$ belongs to $(\F_{01}\setminus\Z)\cup\I$ by \eqref{eq:H} and Lemma~\ref{lem:central}\ref{item:central-incomparable}. Hence, for every $Z\in\Z$ and $Q\in\H_i$, Lemma~\ref{lem:central}\ref{item:central-incomparable} gives $Z\nsubseteq Q$, while Lemma~\ref{lem:central}\ref{item:central-cross} gives $Z\setminus Q\in\D$. Since $\D$ is a downset, all the conditions of Lemma~\ref{lem:dimension} are satisfied. Therefore,
\begin{equation}\label{eq:dimension}
|\D|\ge |\H_i|+\rank M\qquad(i=0,1).
\end{equation}
Summing \eqref{eq:dimension} over $i=0,1$ and using \eqref{eq:Hcount} gives $2|\D|\ge|\F_{01}|+|\I|-|\Z|+2\rank M$. Combining this inequality with \eqref{eq:count}, we obtain
\begin{equation}\label{eq:rank-count}
|\F|\le 2^{n-2}+|\Z|-2\rank M.
\end{equation}
If $\Z=\varnothing$, then $|\F|\le2^{n-2}$. Assume henceforth that $\Z\ne\varnothing$.

For $k\ge2$, \eqref{eq:small} allows us to apply Lemma~\ref{lem:recover} with $\mathcal U=\Z$ and $r=k-2$, yielding
\begin{equation}\label{eq:incidence}
\col W_{k-2}(\Z)\subseteq\col M.
\end{equation}
Each row of $W_{k-2}(\Z)$ is indexed by a $k$-subset of $R$, where $|R|=2k$ and $0\le k-2\le\min(k,|R|-k)$. Lemma~\ref{lem:restricted-rank} and \eqref{eq:incidence} therefore give
\begin{equation}\label{eq:large-rank}
\rank M\ge\rank W_{k-2}(\Z) \ge\frac{\binom{2k}{k-2}}{\binom{2k}k}|\Z| =\frac{k(k-1)}{(k+1)(k+2)}|\Z|.
\end{equation}
For $k\ge6$ the last coefficient is at least $1/2$, so \eqref{eq:rank-count} gives $|\F|\le2^{n-2}$.
\subsection{The cases \texorpdfstring{$n=4,6,8$}{n = 4, 6, 8}}
For $n=4$, we have $|\Z|\le\binom21=2$, so \eqref{eq:rank-count} gives $|\F|\le4+2=6$. For $n=6$, \eqref{eq:small} gives $\varnothing\in\D$. The corresponding column of $M$ is nonzero, since $\Z\ne\varnothing$, and hence $\rank M\ge1$. As $|\Z|\le\binom42=6$, \eqref{eq:rank-count} gives $|\F|\le16+6-2=20=\binom{6}{3}$.

Let $n=8$ and write $z=|\Z|$. If $\Z$ contains no complementary pair, then $z\le\frac12\binom63=10$. Equations~\eqref{eq:rank-count} and \eqref{eq:large-rank} give $ |\F|\le64+z-2\cdot\frac{3z}{10} =64+\frac{2z}{5}\le68<\binom{8}{4}. $ Otherwise, take $Z_0,\overline{Z_0}\in\Z$. If $P\in\H_0\cup\H_1$, applying Lemma~\ref{lem:central}\ref{item:central-cross} to $P$ and each of $Z_0,\overline{Z_0}$ gives
\[
|P|=|P\cap Z_0|+|P\cap\overline{Z_0}|\ge4, \qquad |R\setminus P| =|R\setminus(P\cup Z_0)| +|R\setminus(P\cup\overline{Z_0})| \ge 4,
\]
contrary to $|R|=6$. Hence $\H_0=\H_1=\varnothing$, and \eqref{eq:Hcount} gives $|\F_{01}|+|\I|=z$. Moreover, \eqref{eq:small} applied to the two ordered pairs $(Z_0,\overline{Z_0})$ and $(\overline{Z_0},Z_0)$ shows that $\D$ contains the empty set and all six singletons of $R$. Thus \eqref{eq:count} yields $ |\F|\le64+z-14\le64+\binom63-14=70=\binom{8}{4}$.

\subsection{The cases \texorpdfstring{$n=10,12$}{n = 10, 12}}
In this case, $k=4$ or $5$.

\smallskip\noindent
\emph{Case 1. Every member of $\F_{01}\cup\I$ has size $k$.} Put $\mathcal A=\F_{01}\cup\I$. Suppose first that some $v\in R$ belongs to every member of $\mathcal A$ or to none. Deleting $v$ from the members of $\Z$ gives distinct $(k-1)$-sets or distinct $k$-sets on $2k-1$ points. Their $(k-2)$-inclusion matrix is obtained from $W_{k-2}(\Z)$ by retaining the columns indexed by sets avoiding $v$. Hence \eqref{eq:incidence} and Lemma~\ref{lem:restricted-rank} give
\[
\rank M\ge \rank W_{k-2}(\Z) \ge \frac{\binom{2k-1}{k-2}}{\binom{2k-1}{k-1}}|\Z| =\frac{k-1}{k+1}|\Z|\ge\frac12|\Z|,
\]
and \eqref{eq:rank-count} proves the required bound.

We may therefore assume that every point of $R$ belongs to some but not all members of $\mathcal A$. Put $ \mathcal B=\F_{01}\cap\I$ and $ \mathcal S=\partial^{\le k-2}\mathcal A \cap\partial^{\le k-2}\mathcal A^c.$ By Lemma~\ref{lem:small}, $\mathcal S\subseteq\D$. Since $|\F_{01}|+|\I|=|\mathcal A|+|\mathcal B|$, \eqref{eq:count} shows that it suffices to prove $|\mathcal A|+|\mathcal B|\le2|\mathcal S|$.

For $T\subseteq R$ with $|T|\le k-2$, let $c_T=(\one_{T\subseteq A})_{A\in\mathcal A}$, and let $L$ be the matrix consisting of these columns. Fix $v\in R$ and put $\mathcal S_v=\{T\in\mathcal S:v\notin T\}$. For any $A,A'\in\mathcal A$ and $T\subseteq R$ with $|T|\le k-2$, if $ T\subseteq(A\setminus\{v\})\setminus(A'\setminus\{v\}), $ then $v\notin T$, $T\subseteq A$, and $T\cap A'=\varnothing$, so $T\in\mathcal S_v$.  Since $\mathcal A$ is $k$-uniform, deleting $v$ is injective on $\mathcal A$. We may therefore apply Lemma~\ref{lem:recover} with $V=R\setminus\{v\}$, $r=k-2$, $\mathcal U=\{A\setminus\{v\}:A\in\mathcal A\}$ and $\D=\mathcal S_v$, identifying the rows with those indexed by $\mathcal A$. It follows that the columns $c_T$  of $L$ with $v\notin T$ and $|T|\le k-2$ are spanned by those indexed by $\mathcal S_v$.

Since every member of $\mathcal A$ has size $k$, we also have
\[
c_{T+v}=(k-|T|)c_T -\sum_{j\in R\setminus(T+v)}c_{T+j} \qquad (v\notin T,\ |T|\le k-3).
\]
This expresses each remaining column in terms of columns avoiding $v$, and hence $\rank L\le|\mathcal S_v|$ for every $v\in R$.

Since every point of $R$ belongs to some but not all members of $\mathcal A$, all $2k$ singletons belong to $\mathcal S$. As $\mathcal S$ is a downset, it also contains the empty set. The empty set is counted $2k$ times, each singleton $2k-1$
times, and every other member at most $2k-2$ times.
Hence
\begin{equation}\label{eq:central-upper}
\begin{aligned}
2k\rank L \le\sum_{v\in R}|\mathcal S_v| =\sum_{T\in\mathcal S}(2k-|T|) \le(2k-2)|\mathcal S|+2k+2.
\end{aligned}
\end{equation}
We need one further observation about the size of $\mathcal S$.
By assumption, some member of $\mathcal A$ contains $2k+1$,
and some member avoids $2$. Since $\mathcal A$ is shifted, we have
$\{2,\ldots,k\}\cup\{2k+1\}\in\mathcal A$ and $\{3,\ldots,k+2\}\in\mathcal A$.
Further shifts yield
\[
A_j=\{2,\ldots,k\}\cup\{j\}\in\mathcal A
\qquad(k+3\le j\le2k+1),
\]
and
\[
B_i=\{2,\ldots,k+2\}\setminus\{i\}\in\mathcal A
\qquad(2\le i\le k).
\]
For every such $i,j$, the pair $\{i,j\}$ is contained in $A_j$ and disjoint from $B_i$. Hence these $(k-1)^2$ pairs belong to $\mathcal S$, giving
\begin{equation}\label{eq:central-shadow-size}
|\mathcal S|\ge1+2k+(k-1)^2=k^2+2.
\end{equation}

To bound $\rank L$ from below, we first show that $|X\cap Y|\ge2$ for all $X\in\F_{01}$ and $Y\in\I$. By \eqref{eq:mixed}, $X\cap Y\ne\varnothing$. If $X\cap Y=\{x\}$, then $X-x\subsetneq\overline Y$. But $X-x\in\F_{11}$ by \eqref{eq:shift} and $Y\in\I\subseteq\F_{00}$, contradicting Lemma~\ref{lem:four}(i). Hence $|A\setminus B|\le k-2$ for $A\in\mathcal A$ and $B\in\mathcal B$. For each $B\in\mathcal B$, \eqref{eq:binomial} therefore gives
\[
\sum_{\substack{T\subseteq\overline B\\|T|\le k-2}} (-1)^{|T|}c_T =\bigl(\one_{A=B}\bigr)_{A\in\mathcal A}.
\]
Thus, for every $B\in\mathcal B$, the vector $(\one_{A=B})_{A\in\mathcal A}$ belongs to $\col L$. 
In every linear relation among the rows of $L$, the
coefficient of each row indexed by $\mathcal B$ is zero. Hence
\[
\rank L=|\mathcal B|+\rank L[\mathcal A\setminus\mathcal B].
\] 
 Since $W_{k-2}(\mathcal A\setminus\mathcal B)$ is a submatrix of $L[\mathcal A\setminus\mathcal B]$, Lemma~\ref{lem:restricted-rank} gives
\begin{equation}\label{eq:central-rank}
\rank L\ge|\mathcal B|+\rank W_{k-2}(\mathcal A\setminus\mathcal B) \ge|\mathcal B| +\frac{k(k-1)}{(k+1)(k+2)} (|\mathcal A|-|\mathcal B|) \ge\frac{k}{2(k+1)}(|\mathcal A|+|\mathcal B|).
\end{equation}
The last inequality holds for $k\in\{4,5\}$.

Combining \eqref{eq:central-upper}--\eqref{eq:central-rank}, we obtain
\[
\begin{aligned}
|\mathcal A|+|\mathcal B| \le2|\mathcal S| +\frac{2\bigl((k+1)^2-|\mathcal S|\bigr)}{k^2} \le2|\mathcal S|+\frac{4k-2}{k^2} <2|\mathcal S|+1.
\end{aligned}
\]
Integrality gives $|\mathcal A|+|\mathcal B|\le2|\mathcal S|$, completing Case~1.

\smallskip\noindent
\emph{Case 2. Some member of $\F_{01}\cup\I$ has size different from $k$.} Complementing every member of $\F$ and then reversing the order on $[n]$ preserves shiftedness, the weakly IU property and $|\F|$. On $\F_{01}\cup\I$, these operations replace each set by its complement in $R$ with the order reversed. We may therefore choose $U\in\F_{01}\cup\I$ with $|U|>k$. Fix a $(k-1)$-set $P\subseteq R$ containing $\overline U$. For every $Z\in\Z$, Lemma~\ref{lem:central}(ii) gives $|\overline U\setminus Z|\ge2$, and hence $|P\setminus Z|\ge2$. By \eqref{eq:incidence} and \eqref{eq:rank-count}, it suffices to prove $\rank W_{k-2}(\Z)\ge|\Z|/2$.

Choose $a\in P$ omitted by as many members of $\Z$ as possible. Among the points of $P\setminus\{a\}$, choose $b$ omitted by as many members of $\Z$ containing $a$ as possible. Order the row indices $Z\in\Z$ and the column indices $T\in\binom R{k-2}$ in three groups, placing sets omitting $a$ first, sets containing $a$ but omitting $b$ second, and sets containing both $a$ and $b$ third. Let $B_1,B_2,B_3$ be the corresponding diagonal blocks of $W_{k-2}(\Z)$, and let $m_j$ be the number of rows of $B_j$ for $j=1,2,3$.

Every member of $\Z$ omits at least two points of $P$. Moreover, every member containing $a$ omits at least two points of $P\setminus\{a\}$. The maximal choices of $a$ and $b$ therefore give, respectively,
\begin{equation}\label{eq:omission-counts}
(k-1)m_1\ge2|\Z|, \qquad (k-2)m_2\ge2(m_2+m_3).
\end{equation}

The containment condition $T\subseteq Z$ forces every block above the diagonal of $W_{k-2}(\Z)$ to be zero. Consequently,
\begin{equation}\label{eq:block-rank}
\rank W_{k-2}(\Z) \ge \rank B_1+\rank B_2+\rank B_3.
\end{equation}
The block $B_1$ is the $(k-2)$-inclusion matrix of $\{Z\in\Z:a\notin Z\}$, a family of $m_1$ $k$-sets on the $2k-1$ points of $R\setminus\{a\}$. Deleting $a$ from the row and column indices of $B_2$ gives the $(k-3)$-inclusion matrix of $\{Z\setminus\{a\}:Z\in\Z,\ a\in Z,\ b\notin Z\},$ a family of $m_2$ $(k-1)$-sets on the $2k-2$ points of $R\setminus\{a,b\}$.

Suppose first that $k=4$. Since $|P|=3$ and every member of $\Z$ omits at least two points of $P$, no member of $\Z$ contains both $a$ and $b$. Thus $m_3=0$ and $|\Z|=m_1+m_2$. The bound $3m_1\ge2|\Z|$ from \eqref{eq:omission-counts} therefore implies $m_1\ge2m_2$. Combining \eqref{eq:block-rank} with $m_3=0$ and Lemma~\ref{lem:restricted-rank} applied to $B_1$ and $B_2$ yields
\begin{equation} \label{equ:k4rank}
\rank W_2(\Z) \ge \rank B_1+\rank B_2 \ge \frac{\binom72}{\binom74}m_1 +\frac{\binom61}{\binom63}m_2 \ge\frac12(m_1+m_2) =\frac12|\Z|,
\end{equation}
where the last inequality uses $m_1\ge2m_2$, given by \eqref{eq:omission-counts} with $k=4$ and $m_3=0$.

Now let $k=5$. The three row classes partition $\Z$, so $|\Z|=m_1+m_2+m_3$. The bounds $4m_1\ge2|\Z|$ and $3m_2\ge2(m_2+m_3)$ from \eqref{eq:omission-counts} imply $m_1\ge m_2+m_3$ and $m_2\ge2m_3$.

Every row index $Z$ of $B_3$ contains both $a$ and $b$. Since $|P|=4$ and $Z$ omits at least two points of $P$, we have $Z\cap P=\{a,b\}$. Restrict $B_3$ to the columns indexed by $\{a,b,y\}$, where $y\in R\setminus P$. Deleting $a,b$ from the row and column indices identifies the retained submatrix of $B_3$ with the $1$-inclusion matrix of the $m_3$ triples $Z\setminus\{a,b\}$ on the six points of $R\setminus P$. Lemma~\ref{lem:restricted-rank} therefore gives $\rank B_3\ge\frac{\binom61}{\binom63}m_3$. Combining \eqref{eq:block-rank} with the preceding bound for $\rank B_3$ and Lemma~\ref{lem:restricted-rank} applied to $B_1$ and $B_2$ yields
\begin{equation} \label{equ:k5rank}
\rank W_3(\Z) \ge \frac{\binom93}{\binom95}m_1 +\frac{\binom82}{\binom84}m_2 +\frac{\binom61}{\binom63}m_3 \ge\frac12(m_1+m_2+m_3) =\frac12|\Z|,
\end{equation}
where the last inequality uses $m_1\ge m_2+m_3$ and $m_2\ge2m_3$, both given by \eqref{eq:omission-counts} with $k=5$.

Equations~\eqref{equ:k4rank} and \eqref{equ:k5rank}, together with \eqref{eq:incidence}, give $\rank M\ge|\Z|/2$ for $k=4$ and $k=5$, respectively. Theorem~\ref{thm:main} now follows from \eqref{eq:rank-count}.

\section*{Declaration of AI use}
The proof was obtained using ChatGPT under the guidance of the authors.  The authors then rewrote and restructured the material with the assistance of ChatGPT. The authors independently verified all mathematical statements and proofs and reviewed the entire manuscript and take full responsibility for the correctness and content of the final manuscript.


\begin{thebibliography}{9}
\bibitem{Frankl} P. Frankl, Old and new applications of Katona's circle, \emph{European Journal of Combinatorics} \textbf{95} (2021), 103339.

\bibitem{FT} P. Frankl and N. Tokushige, \emph{Extremal Problems for Finite Sets}, Student Mathematical Library, vol. 86, American Mathematical Society, Providence, RI, 2018.

\bibitem{Gottlieb}
D. H. Gottlieb, A certain class of incidence matrices,
\emph{Proceedings of the American Mathematical Society}
\textbf{17} (1966), 1233--1237.

\bibitem{Gronau} H.-D. O. F. Gronau, On maximal families of subsets of a finite set, \emph{Discrete Mathematics} \textbf{34} (1981), 119--130.

\bibitem{GP} H.-D. O. F. Gronau and Chr. Proske, Maximal families of restricted subsets of a finite set, \emph{Acta Cybernetica} \textbf{5} (1982), 441--451.

\bibitem{Harris} T. E. Harris, A lower bound for the critical probability in a certain percolation process, \emph{Proceedings of the Cambridge Philosophical Society} \textbf{56} (1960), 13--20.

\bibitem{Kleitman} D. J. Kleitman, Families of non-disjoint subsets, \emph{Journal of Combinatorial Theory} \textbf{1} (1966), 153--155.

\bibitem{MS} J. Marica and J. Sch\"onheim, Differences of sets and a problem of Graham, \emph{Canadian Mathematical Bulletin} \textbf{12} (1969), 635--637.
\end{thebibliography}
\end{document}